\documentclass{WMY2000}
\begin{document}
 \input amssym.def
 \input amssym.tex

  \title{L-functions and Random Matrices}
 \author{J. Brian Conrey}
 \institute{American Institute of Mathematics  and Oklahoma State University}
  \maketitle
  
  \section{The GUE Conjecture}
  
\subsection{Introduction}
  In 1972 H.  L.  Montgomery announced a remarkable connection between the
distribution of the zeros of the Riemann zeta-function and the distribution of
eigenvalues of large random Hermitian matrices.  Since then a number of startling
developments have occurred making this connection more profound. 
In particular, random matrix theory has been found to be an extremely useful predictive tool in the theory of L-functions.
  In this article we will try to explain these recent developments and indicate some directions for future investigations.

\subsection{ The Riemann zeta-function}

The
 Riemann zeta-function is defined by
\begin{equation}\zeta(s)=\sum_{n=1}^\infty \frac 1{n^s} \qquad (s=\sigma+it, 
\sigma>1)\thinspace.\end{equation}
It can also be expressed as a product over primes, the Euler product, 
 \begin{equation}\zeta(s)=\prod_p \left(1-\frac 1{p^s}\right)^{-1}\thinspace ,\end{equation}
  for $\sigma>1$.
In 1859 Riemann 
proved that $\zeta(s)$ extends to a meromorphic function on the whole plane with its only 
singularity being a simple pole at $s=1$ with residue 1.
He further proved that it has a functional equation relating the value of $\zeta(s)$
  with the value 
of $\zeta(1-s)$, 
\begin{equation} \zeta(s)=\chi(s)\zeta(1-s)\end{equation}
where $\chi(1-s)=\chi(s)^{-1}=2(2\pi)^{-s}\Gamma(s)\cos(\pi s/2)$. 
He discovered that the distribution of prime numbers is governed by the 
zeros of $\zeta$. He was led to
 conjecture that all of the complex zeros $\rho=\beta+i\gamma$ of 
$\zeta$ have $\beta =1/2$. This assertion is the famous Riemann Hypothesis.
 
We know that $0<\beta<1$ for any complex zero of $\zeta$. Riemann estimated the zero 
counting function 
\begin{equation}N(T)=\#\{\rho=\beta+i\gamma: 0<\gamma \le T\}=\frac{T}{2\pi}\log 
\frac{T}{2\pi e} +O(\log T)\thinspace.\end{equation}
Thus, at a height $T$ the average spacing between zeros is asymptotic to $2\pi/\log T$.
See   \cite{T} for additional background information about $\zeta(s)$.

\subsection{Background}

Montgomery \cite{M1} was studying gaps between zeros of the Riemann zeta-function 
  in an attempt to 
prove that the spacings between consecutive zeros can sometimes be less than 1/2 of the average 
spacing.  
 Such a conclusion would have led to 
a good effective lower bound for the class number of an imaginary quadratic field.
This estimate was not achieved, but through the course of his 
analysis Montgomery  was led to conjecture that
\begin{equation}\lim_{T\to \infty}\frac{1}{N(T)} \sum _{ 0\le \gamma, \gamma' \le T\atop
\frac {2\pi \alpha}{\log T}\le \gamma-\gamma' \le \frac {2\pi \beta}{\log T}
 } 1=\int_\alpha^\beta \left(1-\big(\frac {\sin  \pi x}{\pi x}\big)^2\right)~dx 
\end{equation}
  where $0< \alpha < \beta$ are fixed.

When Montgomery told Freeman Dyson this formula, 
Dyson  responded that the integrand was the pair-correlation function for eigenvalues 
of large random Hermitian matrices, or more specifically the Gaussian Unitary Ensemble, or GUE.

The GUE is the limit as $N\to \infty$ of the probability space consisting of  $N\times N$ Hermitian matrices
$H$  with a probability measure $P(H)dH$ that is invariant under 
conjugation by any unitary matrix $U$. 
 Here $dH=\prod_{j\le k}d\Re H_{jk}\prod_{j<k}d\Im H_{jk}$.
 Mathematical physicists had studied various ensembles 
since the 1950s in connection with work of Wigner in nuclear physics.
Mehta \cite{Me} has given a thorough treatment of the development of the subject.

Montgomery   went on to conjecture that the $n$-correlation function for
zeros of $\zeta$ is the same as that for the GUE; this conjecture came to be known 
as Montgomery's GUE conjecture.

Odlyzko and Sch\"{o}nhage developed an algorithm that allowed for the simultaneous 
calculation of many values of $\zeta(1/2+it)$ for $t$ near $T$ in average time $T^\epsilon$. 
This algorithm allowed Odlyzko \cite{O} to do extensive computations of the zeros of $\zeta$ at a 
height near zero number $10^{20}$; his computations of the pair correlation and nearest 
neighbor spacing for the zeros of $\zeta$ were amazingly close to those for the GUE.  His 
famous pictures added much credibility to Montgomery's conjecture.

\subsection{Further evidence for GUE}

In 1996 Rudnick and Sarnak \cite{RS} made some interesting  progress on the GUE conjecture.
To explain their result, number the ordinates of the zeros of $\zeta(s)$:
$0<\gamma_1\le \gamma_2 \le  \dots $.  Introduce a scaling
$\tilde \gamma =\gamma \frac{\log \gamma}{2\pi}$ so that the $\tilde \gamma$ have asymptotic mean 
spacing 1.
Then Sarnak and Rudnick  proved that
\begin{equation}\lim_{T\to \infty} \frac{1}{T}\sum _{ \gamma_{j_1}, 
\dots ,\gamma_{j_n} \le T\atop
j_m\ne j_n } f(\tilde \gamma_{j_1}, \dots , \tilde \gamma_{j_n})
=   \int_{P_n} W_{U,n}(\vec{x})f(\vec{x}) ~d\vec{x}
 \end{equation}
where
$W_{U,n}(\vec{x})=W_{U,n}(x_1,\dots,x_n)$ is the $n$-correlation function for the GUE (see section 5.5) and where $f$ is any function satisfying
(1) $f(\vec{x}+t(1,\dots, 1))=f(\vec{x})$ for $t\in \Bbb{R}$; (2) $f$ is smooth and symmetric in 
the variables and decays rapidly as $x\to\infty$ in the hyperplane $P_n:=\{(x_1,\dots,x_n):
\sum_{j=1}^n x_j=0\}$; (3) 
the
 Fourier transform $\hat f(\vec{u})$ of $f$ is supported in $\sum_{j=1}^n |u_j|<2$.
The condition (1) assures that $f$ is a function of the differences of the $\gamma_j$.  

Thus, the GUE
conjecture has been proven for all correlation functions, for a limited
class of test functions

\subsection{Other L-functions}

The Riemann zeta-function is the prototype for some extraordinary objects known as 
arithmetic L-functions.  An arithmetic L-function (or L-function for short) has many 
properties in common with the Riemann $\zeta$-function. 
It has a Dirichlet series
\begin{equation}L(s)=\sum_{n=1}^\infty{a_n}{n^{-s}}\thinspace;\end{equation}
it is meromorphic apart from a possible pole at $s=1$; it has a functional equation 
$ \gamma(s)L(s)=\epsilon \overline{\gamma(1-\overline{s})L(1-\overline{s})},$ 
where $\epsilon$, the sign, satisfies $|\epsilon|=1$. $L(s)$  has  an 
Euler product in which the $p$-th factor is the reciprocal of a polynomial in $p^{-s}$.
  Moreover, an arithmetic L-function is 
expected to satisfy the Riemann Hypothesis, that all complex zeros have real part 1/2.  (Note 
that
 our L-functions are normalized so that the 1/2-line is the line of symmetry for the functional 
equation. See \cite{S} for precise definitions.)  An L-function is called primitive if 
it is not the product of two 
L-functions.
Primitive L-functions arise from a variety of contexts: from Dirichlet characters, from 
Mellin transforms of certain cusp forms, from Galois representations, from algebraic 
varieties, etc. However, it is believed that each primitive arithmetic L-function
 is associated to a cuspidal automorphic representation of GL$_m$ over a number field.

Rudnick and Sarnak applied their methods to a fairly general L-function. They proved the analogue of (6)
for any cuspidal 
automorphic L-function over $\Bbb Q$, assuming the Riemann Hypothesis for the L-function,
and with  
obvious changes to reflect 
the appropriate scaling of the zeros. In particular, the answer did not depend in any way on the distribution of the coefficients of the 
particular L-function.

The GUE Conjecture is now seen as
a universal law governing the distribution of zero spacings for
all arithmetic L-functions.

\section{Families}
The realization that the compact classical groups of 
 matrices play a role in the theory of L-functions arose 
 through seminal work of Katz and Sarnak \cite{KS1}, \cite{KS2}.

\subsection{Function Field Analogues}
Katz and Sarnak investigated the distribution of zeros of function field 
zeta-functions. Consider a curve $f(x,y)=0$ where $f$ is a 
  polynomial with integer coefficients.  
The zeta-function for $f$ over a finite field $\Bbb{F}_q$ can be obtained in a simple way from 
the 
generating function of the numbers $N_n$ of points on that curve in the finite field extensions
$\Bbb{F}_{q^n}, n=1,2,\dots.$  It is known that this zeta-function is a rational function 
whose numerator is a polynomial with 
integer coefficients, degree 2$g$ where $g$ is the genus of the curve, and 
that it satisfies 
the Riemann Hypothesis that all zeros have modulus $1/\sqrt{q}$.  One can order these zeros on 
the circle $|z|=\sqrt{q}$ in terms of their angles measured from the positive axis. One
  then considers statistics of the angles.  

Katz and Sarnak proved that, after proper normalization,
the $n$-correlation function of these angles, as $g$ and $q$ tend to infinity, is exactly 
the $n$-correlation of GUE.
They also proved that the nearest neighbor, or consecutive spacing,  statistic for these
 angles
is the same as that for GUE.

The method of proof involved working with subgroups of U($N$), the group of 
$N\times N$ unitary matrices   with 
the Haar measure (see sections 5.1 and 5.2). It 
had been  shown by Dyson that the $n$-correlation and the spacing statistics are the 
same for GUE as for the appropriately scaled limit as $N\to \infty$ for    $U(N)$ 
with Haar measure. (In the physics literature the limit of $U(N)$ is called the Circular Unitary 
Ensemble or CUE; see \cite{Me}.)

The function field zeta-functions that Katz and Sarnak were working with were characteristic 
polynomials of   matrices from subgroups of U($N$). To deduce their result they applied a theorem of Deligne about the 
equidistribution of function field zeta-functions among the characteristic polynomials of 
conjugacy classes in these subgroups.

\subsection{Introduction of Families}
Katz and Sarnak investigated the robustness of their theorem.  They conjectured that the conclusion
  remains true  under the weaker hypothesis that $q$ is held fixed and $g\to \infty$.  
  On the other hand, they provided
examples of sequences of curves of increasing 
genus for which the spacing statistic is   not the GUE consecutive spacing statistic.  

They also considered special families such as (a) curves of the form $y^2=f(x)$ with 
all squarefree, monic $f$, 
and (b) quadratic twists of curves of the form $y^2=x(x-1)(x-t)$; that is
$\Delta y^2=x(x-1)(x-t)$ where $\Delta$ is a squarefree, monic polynomial over ${\Bbb F}_q$. 
 When $q$ and the degree of $\Delta\to \infty$ they discovered that 
 the $n$-correlation and spacing statistics 
again matched up with CUE.
However, the reason  in each case is different depending on the geometric monodromy group.
In case (a) they computed 
that   the monodromy group was the
  symplectic group and in case (b)  the orthogonal group. 

We denote by USp($N$) (if $N$ is even) the unitary symplectic matrices, 
and by O($N$) the orthogonal matrices. 
These are subgroups of U($N$) and are equipped with their own Haar
 measure.  
 
 Katz and Sarnak proved that the $n$-correlation statistic and the spacing 
statistic for the limits of U($N$), USp($N$), and O($N$) are all the same as that for GUE.

By contrast, the Circular Orthogonal Ensemble (COE) and the Circular Symplectic Ensemble (CSE),
which are well-known ensembles in mathematical physics 
 (see \cite{Me}) have the same underlying symmetry groups as O and Sp,
 but have a measure different from the Haar 
measure.  They have different $n$-correlation and spacing statistics than those 
of GUE and CUE.

From now on we will use the letters O, Sp, U when referring to a statistic associated with the 
 limits of  O($N$), USp($N$), and U($N$). 
These statistics are computed for $N\times N$ matrices of 
the appropriate subgroup of U$(N)$ with their respective Haar 
measure and then $N\to \infty$ with an appropriate scaling limit in 
order to determine the statistics (see 5.4, 5.5, and 5.6).

Since orthogonal and symplectic matrices have eigenvalues which occur in complex conjugate 
pairs, it is clear that the eigenvalue 1 plays a special role for these matrices, whereas it 
does not for unitary matrices.
Indeed, there are statistics which   differentiate   O, Sp, and U.  In particular, 
the eigenvalue nearest to 1 is such an example. It turns out that this statistic is 
different for all three symmetry types.  More generally the $j$-th eigenvalue nearest to 1 is a 
statistic that is dependent on the group (for a summary of the four statistics 
of interest to us, see section 
5.3).

 Another statistic that is different for all three symmetry types is ``level-density''.  The 
$n$-level density function is
 obtained from  summing a test function at $n$-distinct eigenvalues.
 For comparison purposes, note that the $n$-correlation function is obtained by summing a test function of $n$-variables which 
is a function depending only on the differences of the arguments at $n$ distinct eigenvalues.
 See sections 5.5 and 5.6 for more on these statistics.
 
Katz and Sarnak showed that for their special families of function field zeta-functions 
these new statistics ($j$-th nearest eigenvalue to 1 and $n$-level density) match with the 
appropriate statistics from Sp and O, which they had computed. 

In general, they found that if they could compute the geometric monodromy group associated 
with a family of function field zeta-functions, and if that monodromy group was U, Sp or O, then 
all four of the statistics we have been discussing 
for the function field zeta-functions could be proven to match with the appropriate 
statistic from U, Sp or O. 

Thus, we say that the symmetry type for the family   of 
curves of the form $y^2=f(x)$ the symmetry type is Sp; and for the family of quadratic twists of 
$y^2=x(x-1)(x-t)$ 
the symmetry type is O.   

Katz and Sarnak also have examples of families of function field zeta-functions where the symmetry type is U.

\subsection{L-Functions Over Number Fields}
Katz and Sarnak speculated that their results for function field zeta-functions would have 
implications for L-functions over number fields. Two collections  of L-functions which present 
themselves as natural analogues to the families above are 
(a$'$) the collection of all Dirichlet L-functions $L(s,\chi_d)$ where $\chi_d(n) $ is a real primitive Dirichlet character with conductor $|d|$ (so that $d$ runs 
through the set of fundamental discriminants of quadratic number fields) and (b$'$)  the 
collection of L-functions associated with primitive Hecke newforms of a fixed weight. The 
analogy with the function field zeta-functions is as follows.  The zeta-function for a 
member of the Sp family (a) above is obtained by counting solutions to the equation $y^2=f(x)$ 
over a finite field. These solutions are counted in a field with $p$ elements by the sum 
$\sum_{a=1}^p \chi_p(f(a))$ where the real Dirichlet character $\chi_p$ modulo $p$ 
appears.  This suggests that the family (a$'$) could have a symmetry type Sp. Similarly, the 
zeta-function of a member of the O family (b) above is the reciprocal of the $p$-th
  factor in the Euler product for the L-function 
of the elliptic curve defined through twists of the equation $y^2=x(x-1)(x-t)$, the 
Legendre family of elliptic curves. These L-functions are known (by the solution of the 
Taniyama - Shimura conjecture) to be associated to primitive Hecke newforms of weight 2. 
Thus, the collection (b$'$) could well be a family with symmetry type O.

The first evidence that (a$'$) is a family with symmetry type Sp  came from Michael Rubinstein's 
thesis \cite{R}. He computed the lowest lying zero of $L(s,\chi_d)$ and 
the data matched well with the eigenvalue nearest 1 for symplectic matrices. He examined 
theoretically the $n$-level density of zeros and showed
 (for test functions with 
restrictions on the support of their Fourier transforms) that the $n$-level density functions 
for the zeros of $L(s,\chi_d)$ are identical with the $n$-level density functions for Sp. 
Precisely, 
we index the ordinates of the zeros of $L(s,\chi_d)$
as
\begin{equation}0\le \gamma_1^{(d)}\le \gamma_2^{(d)}\le...\end{equation}
and scale using
\begin{equation}\tilde \gamma^{(d)}=\gamma^{(d)} \frac{\log \gamma^{(d)}}{2 
\pi}\thinspace.\end{equation}
Let $D^*=\sum_{|d|\le D} 1$. Then Rubinstein proved that   
\begin{equation} \frac{1}{D^*}\sum_{|d|\le D} \sum _{ \gamma_{j_1}, \dots \gamma_{j_n} 
 \atop
  j_m\ne j_n } f\big(\tilde \gamma_{j_1}^{(d)}, \dots , \tilde \gamma_{j_n}^{(d)}\big)
\to   \int_{{\Bbb R}^n} W_{Sp,n}(\vec{x})f(\vec{x})~d\vec{x}\thinspace,
 \end{equation}
as $D\to \infty$ where $W_{Sp,n}$ is as in section 5.5,
provided that $f$ is a Schwarz function such that  the support of $\hat f(u)$ is contained in
$\sum_{j=1}^n |u_j|<1$.
$W_{Sp,n}$ is called the $n$-level density function for the symplectic group.

Rubinstein also found evidence that another family has a symmetry type O. To
describe this family let
$\Delta(z)=\sum_{n=1}^\infty \tau(n)\mbox{e}(nz)$ where $\tau$ is Ramanujan's 
tau-function.  It is well-known that $\Delta$ is a primitive Hecke newform of weight 12 and 
level 1. The family (c$'$) is obtained from quadratic twists of the L-function associated with 
$\Delta$ namely
\begin{equation}L(\Delta, s,\chi_d)=\sum_{n=1}^\infty 
\frac{\tau(n)n^{-11/2}\chi_d(n)}{n^s}\thinspace .\end{equation}
If $d<0$ then this L-function automatically vanishes at $s=1/2$ because the associated 
functional equation occurs with $\epsilon=-1$. Rubinstein computed the lowest   zero of 
$L(\Delta,s,\chi_d)$ for $d>0$ and the lowest   zero above the real axis for $L(\Delta, 
s, \chi_d)$ with $d<0$. At this point, we should mention that an orthogonal matrix has 
determinant +1 or $-1$. So there are actually two symmetry types O$^+$ and O$^-$ (and the statistics of O are an average of the statistics of these).
  Katz and Sarnak had computed the statistics 
(neighbor spacing, correlations, density, eigenvalue nearest 1) for these two symmetry types 
as well.
Rubinstein found that the lowest lying zero of L-functions in (c$'$) with $d>0$ followed O$^+$ 
while the lowest lying zero above the real axis for L-functions with $d<0$ followed O$^-$. 

Rubinstein considered more generally the twisting of an arbitrary arithmetic L-function 
$L(f,s)$ associated with an automorphic form $f$ on $GL_m$ by quadratic characters $\chi_d$.  
The results here divide into three cases which have to do with signs of the functional 
equations. 
If
the L-function $L(f,s,\mbox{sym}^2)$ associated with the symmetric square 
of $f$ is entire, the functional equations for $L(f,s,\chi_d)$ will always have $\epsilon=+1$
and the average over $d$ is exactly as above with 
the $n$-level density function being $W_{Sp,n}$.
In the case  that the symmetric square $L(f,s,\mbox{sym}^2)$ has a pole
(at $s=1$) then the sign of the functional equation for $L(f,s,\chi_d)$ is +1 for even 
characters $\chi_d$ and $-1$ for odd characters $\chi_d$.
Rubinstein averages over these two cases separately and discovers 
that the first case yields an $n$-level density function
$W_{O^+,n}$  
and in the second case 
an $n$-level density function $W_{O^-,n} $ (see section 5.5). 
In each of these cases a restriction is placed on the support of the Fourier transform of 
$f$.

\subsection{The Diagonal Terms}
Rubinstein's work gave impressive confirmation of the theory, but the severe restriction on 
the support of the Fourier transform is worth investigating. In fact, this restriction 
occurs right at the place where something interesting is happening with the Fourier transfom 
of the density function. Returning for a moment to the case of the Riemann zeta-function, we 
can illustrate this idea. 

Montgomery's original theorem involved 
\begin{equation}F(\alpha, T)=\frac{1}{N(T)}\sum_{0<\gamma,\gamma' \le T} 
T^{i\alpha(\gamma-\gamma')}w(\gamma-\gamma')\thinspace.\end{equation}
Here $w$ is a weight function that concentrates at the origin. (Montgomery used 
$w(u)=4/(4+u^2)$.)
Assuming the Riemann Hypothesis he showed that 
\begin{equation} F(\alpha, T)= T^{-2\alpha}(1+o(1))+|\alpha| +o(1)\end{equation}
uniformly for $-1+\epsilon <\alpha < 1-\epsilon$ for any $\epsilon >0$. 
 
By the definition of $F$,
\begin{equation}\frac 1{N(T)} \sum_{0<\gamma,\gamma'\le T} r(\gamma-\gamma') 
w(\gamma-\gamma')=\int_{-\infty}^\infty \hat{r}(\alpha)F(\alpha,T) ~d\alpha 
\thinspace.\end{equation}
Thus, Montgomery's theorem gives information about the average behavior of differences 
between the zeros
for test functions $r$ with the support of $\hat{r}$ contained in $(-1,1)$.
Montgomery went on to conjecture (based on considerations of the behavior of prime pairs) 
that for $|\alpha|>1$ one has $F(\alpha,T)=1+o(1)$; this assertion implies (1). 
 Thus, $F$ is not differentiable  at $\alpha=1$. 
In the proof of Montgomery's theorem (via the explicit formula and
 the mean value theorem (23) for Dirichlet polynomials) for $|\alpha|<1$, the main term of $F$
   arises from the ``diagonal'' 
contributions of the mean square of a Dirichlet polynomial (i.e. the terms $m=n$ in the 
integral $\int_0^T a_m \overline{a_n}(m/n)^{it}~dt$).  For $\alpha >1$, the off-diagonal 
terms (i.e. $m\ne n$) contribute to the main-term. 

A similar situation arises in the work of Rudnick and Sarnak and in the work of Rubinstein.  
All of the proofs of these theorems are valid only
 in the range where the diagonal terms dominate.

\subsection{Beyond the Diagonal}
So far we have seen that the theory of families is confirmed by numerical data as well as 
theoretical data up to the diagonal. 
Bogolmony and Keating gave a heuristic derivation of all of the GUE conjecture (i.e. all the 
$n$-level correlations) based 
on  Hardy-Littlewood type conjectures for  pairs of primes and pairs of almost primes;
 this work shows how   the 
off-diagonal terms 
potentially contribute. 

\"{O}zl\"{u}k \cite{Oz} proved an analogue for all primitive Dirichlet characters
for the  pair correlation theorem (12); he obtained a result for $|\alpha|<2$.  This was
the first example of going beyond the diagonal.  See also \cite{OS} and \cite {IS2}.

It is of great interest that Iwaniec, Luo, and Sarnak \cite{ILS} have succeeded
 in going beyond the 
diagonal in several examples which represent   three of the symmetry types. 
They work with the 1-level density functions assuming only that the Riemann Hypothesis holds for all of the 
L-functions in question. The 1-level density functions may be obtained from 
  $W_n$ by taking $n=1$.  They are
\begin{equation}W(\mbox{O})(x)=1+\frac 12 \delta_0(x)\thinspace,\end{equation}
\begin{equation}W(\mbox{O}^+)(x)=1+\frac{\sin2\pi x}{2 \pi x}\thinspace,\end{equation}
\begin{equation}W(\mbox{O}^-)(x)=1-\frac{\sin2\pi x}{2 \pi 
x}+\delta_0(x)\thinspace,\end{equation}
\begin{equation}W(\mbox{Sp})(x)=1-\frac{\sin2\pi x}{2 \pi x}\thinspace,\end{equation}
\begin{equation}W(\mbox{U})(x)=1\thinspace.\end{equation}
The families considered by Iwaniec, Luo, and Sarnak are related to modular forms. Let
$H_k^*(N)$ denote the set of holomorphic newforms $f$ of weight $k$ and level $N$. Let 
$H_k^+(N)$ denote the weight $k$ level $N$ newforms for which the associated L-function has 
a + in its functional equation, and $H_k^-(N)$ is the subset of $f$ for which $L(f,s)$ has a 
$-$ in its functional equation.  
Let $M^*(K,N)$ be the union of the $H^*(k,N)$ for $k\le K$ and similarly define $M^+$ and 
$M^-$. 
They consider the low lying zeros of $L(f,s)$ as $f$ varies through one of the $M$-sets. The 
average spacing for all the zeros of all the $L(f,s)$ with $f\in H^*(k,N)$ up to a fixed 
height $t_0$ is
asymptotic to $2\pi/\log(k^2N)$. Let 
$\phi$ be a test function which is even and rapidly decaying. 
They proved that if the support of $\hat{\phi}$ is contained in $(-2,2)$, then
\begin{equation}\frac{1}{|M^*(K,N)|}\sum_{f\in 
M^*(K,N) \atop \gamma_f} \phi\left(\frac{\gamma_f \log 
k^2N}{2\pi}\right)\to \int_{-\infty}^\infty\phi(x)W(\mbox{O})(x)~dx\thinspace.\end{equation}
as $KN\to\infty$.  Similar statements hold with $M^*$ replaced by $M^+$ and $M^-$ and O replaced 
by O$^+$ and 
O$^-$.
 
It should be pointed out that the Fourier transforms of the  density functions 
$W(\mbox{O})(x)$,
$W(\mbox{O}^+)(x)$, and $W(\mbox{O}^-)(x)$ all agree in the diagonal range; so it is only 
when one goes beyond the diagonal that the distinguishing features of these three symmetry types 
becomes apparent.

Iwaniec, Luo, and Sarnak also consider the symmetric square L-functions of the $f\in M^*$ 
and verify that the above statements hold with symmetry type Sp and the support  of $\hat{\phi}$ in
(-3/2,3/2). Also, the average zero spacing is $2\pi/\log(k^2N^2)$ so $N$ should be replaced 
by $N^2$ in the argument of $\phi$ in the left hand side of (20). 

Thus, the theoretical and numerical evidence 
that the zeros of families of L-functions depend on the symmetry
type of the family is pretty convincing.

\section{Moments }
So far we have seen that the eigenvalues of matrices from unitary groups are excellent 
models for zeros of families of L-functions.  Now we want to take the matrix models a 
significant step further and argue that the characteristic polynomials of these matrices on 
average reveal very important features of the value distribution of the L-functions in the 
family.

\subsection{Moments of the Riemann zeta-function}
We first look at the situation of the Riemann zeta-function and its moments.

To give some background, we cite the theorem of Hardy and Littlewood:
\begin{equation}\frac{1}{T}\int_0^T|\zeta({\textstyle \frac 12}+it)|^2~dt \sim \log 
T \end{equation}
   and the theorem  of Ingham:
\begin{equation}\frac{1}{T}\int_0^T|\zeta({\textstyle \frac 12}+it)|^4~dt \sim \frac {1}{2\pi 
^2}\log^4 
T\thinspace.\end{equation}
The asymptotics of no other moments (apart from the trivial 0-th moment) are known.
In general it has been conjectured that
\begin{equation} \frac{1}{T}\int_0^T|\zeta({\scriptstyle \frac  12}+it)|^{2k}~dt \sim c_k 
\log^{k^2} T\thinspace.\end{equation}

The basic tools for investigating mean-values in $t$-aspect are the mean value theorem for
Dirichlet polynomials (due to Montgomery and Vaughan):  \begin{equation} \int_0^T\bigg|\sum_{n\le N}a_n
n^{it}\bigg|^2~dt=\sum_{n\le N} (T+O(n))|a_n|^2 \end{equation} and some sort of formula
expressing the function in question in terms of Dirichlet polynomials (such as  an approximate
functional equation) such as \begin{equation} \zeta(s)^k=\sum_{n\le \tau^k}\frac{d_k(n)}{n^s}+\chi(s)^k
\sum_{n\le \tau^k}\frac{d_k(n)}{n^{1-s}} +E(s)\thinspace,\end{equation} where $E(s)$ should be small on average,
$s=1/2+it$, $\tau=\sqrt{t/(2\pi)}$, $\chi(s)$ is the factor from the functional equation (3), and
where \begin{equation} \zeta(s)^k=\sum_{n=1}^\infty \frac{d_k(n)}{n^s}\qquad (\sigma>1)
\end{equation} so that $\zeta(s)^k$ is the generating function for $d_k(n)$.  Note that the
mean-value theorem for Dirichlet polynomials detects only diagonal contributions.

 Conrey and Ghosh \cite{CG2} gave the moment
conjecture a more precise form, namely that there should be a factorization 
\begin{equation}c_k
=\frac{g_k a_k}{\Gamma(1+k^2)}\end{equation} 
where  
\begin{equation}a_k=\prod_p\left(1-{  \frac 1p}\right)^{k^2} 
\sum_{j=0}^\infty \frac{d_k(p^j)^2}{p^j} \end{equation} is an arithmetic factor and $g_k$, a
geometric factor, should be an integer. 
 Note that by the mean-value theorem for Dirichlet
polynomials it is not difficult to show that 
\begin{equation} \frac 1T \int_0^T\bigg|\sum_{n\le
x}\frac{d_k(n)}{N^{1/2+it}}\bigg|^2~dt \sim \frac{a_k (\log x)^{k^2}}{\Gamma(1+k^2)}
\thinspace,\end{equation}
provided that $x=o(T)$. 
Thus,  an interpretation of $g_k$
is \begin{equation}g_k=\lim_{T\to \infty}\frac{\int_0^T|\zeta^k({\scriptstyle \frac
12}+it)|^{2}~dt} {\int_0^T\left|\sum_{n\le T}\frac {d_k(n)}{n^{ 1/2 +it}} \right|^2~dt}
\end{equation} assuming that the limit exists, so that $g_k$ represents the `number' of Dirichlet
polynomial approximations to $\zeta(s)^k$ of length $T$ required to measure the  mean square
of $\zeta(s)^k$.  

 In
this notation, the result of Hardy and Littlewood is that $g_1=1$ and Ingham's result is that
$g_2=2$.  

These results can be obtained essentially from the mean-value theorem for Dirichlet polynomials.
To go beyond the fourth moment requires taking into account off-diagonal contributions.
Goldston and Gonek \cite{GG} describe a precise way to transform information about
coefficient correlations $\sum_{n\le x} a(n)a(n+r)$ into a formula for 
the mean square of a long Dirichlet polynomial $\sum_{n\le x} a(n)n^{-s}$ where $x$ is bigger than 
the length of integration.
 
Using Dirichlet polynomial techniques Conrey and Ghosh \cite{CG1}
conjectured that $g_3=42$ and Conrey and Gonek \cite{CGo} conjectured 
that
$g_4=24024$.
 Meanwhile, Keating and Snaith \cite{KeSn1} computed the
moments of characteristic polynomials of matrices in $U(N)$ and 
found that for any real $x$ and any complex number $s$,
\begin{equation}
M_{U,N}(s)=\int_{U(N)} |\det(A-I\exp(-i x))|^{2s} 
~dA=\prod_{j=1}^N\frac{\Gamma(j)\Gamma(j+2s)}{\Gamma(j+s )^2}\thinspace,
\end{equation}
where $dA$ denotes the Haar measure for the group $U(N)$ of $N\times N$ unitary matrices. 
 To do this calculation, they made use of Weyl's formula for the Haar measure 
 (see section 5.3) and Selberg's integral (see section 5.6).
They also showed that 
\begin{equation} \lim_{N\to \infty} \frac{M_N(s)}{N^{s^2}}=\frac{G(1+s)^2}{G(1+2s)}
\thinspace,
\end{equation}
where $G(s)$ is Barnes' double Gamma-function which satisfies $G(1)=1$ and $G(z+1)=\Gamma(z)G(z)$.
Note that for $s=k$ an integer,
\begin{equation}
\frac{G(1+k)^2}{G(1+2k)}=\prod_{j=0}^{k-1}\frac {j!}{(j+n)!}
\thinspace.
\end{equation}

For $k=1,2,3$ the above is 1/1!,2/4!, 42/9!
 in agreement with the theorems of Hardy and Littlewood, and Ingham and
 the conjecture of Conrey and Ghosh.  Keating and Snaith argued that
one should thus model the moments of the zeta-function from 0 to $T$ 
by moments of characteristic 
polynomials of unitary matrices of size $N\sim \log T$. (More precisely, one should take 
$N$ to be the integer nearest to $\log \frac{T}{2\pi}$).
They then conjectured 
that \begin{equation}g_k=k^2!\prod_{j=0}^{k-1} \frac{j!}{(j+k)!}\end{equation} for integer $k$.
   The initial public
announcements   of the conjectures of Conrey and Gonek (that $g_4=24024$) and of Keating and 
Snaith (  $g_k$   for all real $k\ge-1/2$) occurred at 
the Vienna conference on the Riemann Hypothesis only 
moments after it was checked that the Keating and 
Snaith conjecture does indeed predict that $g_4=24024$.

\subsection{Moments of L-functions at 1/2}

Subsequently, Conrey and Farmer \cite{CF} analyzed known results for moments of L-functions
at 1/2 and made a general conjecture. 
 (These moments had been considered by a number of authors; see especially \cite{GV}.)
  The conjecture has the  shape
\begin{equation}\frac{1}{X^*}\sum_{ f\in {\mathcal F}\atop c(f)\le X} L(f,1/2)^k\sim
 \frac{g_k a_k}{\Gamma(1+q(k))} (\log X)^{q(k)} \end{equation}
for some $a_k$, $g_k$, and $q(k)$ where ${\mathcal F}$ is a family of $f$ parametrized by 
the conductor $c(f)$, and $X^*=\sum_{c(f)\le X}1$. 
The observations of Conrey and Farmer were that 
$g_k$ and $q(k)$ depend only on the symmetry type of the family, and that $a_k$ depends 
on the family itself, but is explicitly computable in any specific case.
Thus, the conjecture is that $q(k)=k^2$ for a unitary family, $q(k)=k(k+1)/2$ for a 
symplectic family, and
$q(k)=k(k-1)/2$ for an orthogonal family.
The values of $g_k$ were left unspecified, but were then predicted as before from random matrix theory
 by Keating and Snaith 
\cite{KeSn2}
and independently by Brezin and Hikami \cite{BH} by computing moments of characteristic polynomials of 
matrices from O$(N)$ and from USp(2$N)$. Each is a quotient
 of products of Gamma-functions.

Thus, for the family of Dirichlet L-functions $L(s,\chi_d)$ with a real primitive Dirichlet  
character $\chi_d$ modulo $d$ we have the following results ($D^*=\sum_{|d|\le D}1$):
 Jutila \cite{J} proved that
\begin{equation} \frac{1}{D^*}{\sum_{|d|\le D}}  L({\scriptstyle \frac  12} ,\chi_d) \sim   
a_1\log( D^{\scriptstyle \frac  12 })\end{equation}
and 
\begin{equation} \frac{1}{D^*}{\sum_{|d|\le D}}  L^2({\scriptstyle \frac  12} ,\chi_d) \sim 
2\frac{ 
a_2 
\log^3(D^{\scriptstyle \frac  12} )}{3!}\thinspace.\end{equation}
Soundararajan \cite{So1} showed that
\begin{equation} \frac{1}{D^*}{\sum_{|d|\le D}}  L^3({\scriptstyle \frac  12} ,\chi_d) \sim 
16
 \frac{a_3\log^6 (D^{\scriptstyle \frac  12} )}{6!}\end{equation}
and conjectured that
\begin{equation} \frac{1}{D^*}{\sum_{|d|\le D}} L^4({\scriptstyle \frac  12} ,\chi_d) \sim 
768 
\frac{ a_4
\log^{10}(D^{\scriptstyle \frac 12} )}{10!}\thinspace.\end{equation}

The general   conjecture   coming from random matrix theory (see section 5.7), which agrees with the above, is:
\begin{equation} \frac{1}{D^*}{\sum_{|d|\le D}} L^k({\scriptstyle \frac  12} ,\chi_d)
 \sim  \prod_{\ell=1}^{k}\frac{\ell!}{2\ell!}
a_k\log^{k(k+1)/2}(D)\end{equation}
where 
\begin{equation}a_k=\prod_p\frac{\left(1-{\scriptstyle \frac  1p}\right)^{{\scriptstyle \frac 
{k(k+1)}{2}}}}
{\left(1+{\scriptstyle \frac  1p}\right)}\left(\frac{\left(1-{\scriptstyle \frac 
{1}{\sqrt{p}}}\right)^{-k}+
\left(1+{\scriptstyle \frac {1}{\sqrt{p}}}\right)^{-k}}{2}+\frac 1p\right)\thinspace.\end{equation}

 An example of an orthogonal family where several moments are known arises from  $\mathcal F_q$  the set of 
 primitive cusp forms of weight $2$ and level $q$
($q$ prime). Then,
from results of Duke \cite{D}, Duke, Friedlander, and Iwaniec \cite{DFI}, Iwaniec and 
Sarnak \cite{IS1}, and Kowalski, Michel, and VanderKam \cite{KMV1} and \cite{KMV2}, we have
\begin{equation}\frac{1}{|\mathcal F_q|} \sum_{f\in \mathcal F_q} L(1/2,f)\sim a_1\end{equation}
\begin{equation}\frac{1}{|\mathcal F_q|} \sum_{f\in \mathcal F_q} L^2(1/2,f) \sim 2 a_2 \log 
q^{\scriptstyle \frac  12}\end{equation}
 \begin{equation}\frac{1}{|\mathcal F_q|} \sum_{f\in \mathcal F_q} L^3(1/2,f)  \sim 8a_3 \frac{
  \log^3 q^{\scriptstyle \frac 12}}{3!} \end{equation}
\begin{equation}\frac{1}{|\mathcal F_q|} \sum_{f\in \mathcal F_q} L^4(1/2,f)\sim   
128 a_4 \frac{\log^6 q^{\scriptstyle \frac  12}}{6!}\end{equation}
where
$a_1=\zeta(2),$
\begin{equation}a_2= \zeta(2)^2 \prod_p \left(1+{\scriptstyle \frac {1}{p^2}}\right)
 \end{equation}
 \begin{equation}a_3=\zeta(2)^3\prod_p\left(1-{\scriptstyle \frac  1p}\right)\left(1+{\scriptstyle 
\frac  
1 p }
 +{\scriptstyle \frac  4 {p^2}}+{\scriptstyle \frac  1 {p^3}}+{\scriptstyle \frac  1 
{p^4}}\right)
 \end{equation}
 \begin{equation} a_4 = \zeta(2)^5\prod_p\left(1-{\scriptstyle \frac {1}{p}}\right)^3
 \left(1+{\scriptstyle \frac  3 p} +{\scriptstyle \frac  {11}{p^2}}+
 {\scriptstyle \frac  {10}{p^3}}+{\scriptstyle \frac {11}{p^4}}+
 {\scriptstyle \frac {3}{p^5}}+{\scriptstyle \frac {1}{p^6}}\right)\thinspace \end{equation}
We have not found a simple expression for $a_k$, though it 
can be determined explicitly for each $k$.
 As before, a general conjecture is:
\begin{equation}\frac{1}{|\mathcal F_q|}
 \sum_{f\in \mathcal F_q} L^k(1/2,f)  \sim  2^{k-1} \prod_{\ell=1}^{k-1}\frac{\ell!}{2\ell!}
  a_k \log^{k(k-1)/2} q\thinspace.
   \end{equation}
 
We believe that in  formulas for moments of L-functions over a family
the power of the 
log of the conductor and the value of $g_k$ should only depend on 
the symmetry type of the family and that the value of $a_k$ will 
depend on  the family but can always be determined explicitly.
 
\section{Further directions}
 In this section we mention some questions where further research is desirable.
 \subsection{ Full moment conjecture}  What are the lower order terms in the moment formulae for 
$|\zeta(1/2+it)|^{2k}$ and for $L(1/2)^k$?  These are known in a few instances (see \cite{In}, \cite{C} for the second and fourth moments of $\zeta(s)$) but not in general. 
The difficulty is that random matrix theory does not ``see'' the contribution of the arithmetic factor 
$a_k$.  Lower order terms will likely involve a mix of derivatives of $a_k$ and secondary terms
 from the moments of the characteristic polynomials of matrices.  In general, a better 
understanding of how $\zeta(s)$ is modeled by a characteristic polynomial of a certain type of matrix 
is needed; how do the 
primes come into play?  Perhaps we should think of $\zeta(1/2+it)$ as a partial 
Hadamard product over zeros multiplied by a partial Euler product.  Perhaps these two parts behave 
independently, and the Hadamard product part can be modeled by random matrix theory.

\subsection{ Distribution of Values}  Keating and Snaith 
compute  explicit formulas for the $s$-th  moment of the characteristic
 polynomials of matrices from O($N$), Sp($N$), and U($N$).  Consequently  
  the value distributions for these characteristic 
polynomials can be explicitly computed;
 they involve the Fourier transform of the $s$-th moment.  Preliminary investigations
  indicate that there is  
 a good fit between 
 the random matrix formulae and numerical data. One particularly interesting 
feature of this investigation involves the understanding of zeros which occur exactly at 1/2. 
 These seem to occur only for L-functions in an orthogonal family.  If the sign of the functional 
equation is $-1$, there is automatically a zero at 1/2. The interesting situation 
is when the sign 
is + and there is still a zero; for example if E is an elliptic curve defined over $\Bbb{Q}$
 then examination of the distribution of values of O$^+$ suggest that twists $L(E,s,\chi_p)$ of 
the L-function by quadratic characters seem to vanish for about
 $X^{3/4}(\log X)^{-5/8}$ values of $|p|<X$ with sign +1 in the functional equation. 
 
\subsection{Extreme Values} 

How large is the maximum value of $|\zeta(1/2+it)|$ for $T<t<2T$?
 It is known that the Riemann Hypothesis implies that the maximum is at most
 $\exp(c\log T/ \log\log T)$ for some $c>0$. It is also known that the 
 maximum gets as big as
 $\exp(c_1(\log T/\log \log T)^{1/2})$ for a sequence of $T\to \infty$ for some $c_1>0$.
 It has been conjectured that the smaller bound
  (the one that is known to occur) is closer to the 
truth.
 However, the new conjectures about moments suggest that it may be the larger.
 
 This question has a number of equivalent and analogous (for L-values) formulations. How big can 
$S(T):=\frac1\pi \arg \zeta(1/2+iT)$ be? Assuming the Riemann Hypothesis, it is known that
 $S(T)\ll \log T/\log \log T$ but that infinitely often it is bigger than
 $c(\log T/\log \log T)^{1/2}$ for some $c>0$ Which is closer to the truth?
    What is the maximum size of the class number of an imaginary 
    quadratic field (as a function of 
the discriminant) (see \cite{Sh}for a discussion and numerical investigation.)? How big can 
 the least quadratic
 non-residue of a given prime $p$ be ($\log p$ 
or $(\log p)^{2-\epsilon}$? See \cite{M2} for a discussion of this question?
 What is the maximal order of vanishing of an L-function 
 at 1/2? In terms of the conductor $N$, 
can it be as big as $\log N/\log \log N$
 or is it at most the square root of that, or something entirely different?
All of these questions are related, at least by analogy, and   
 they may all have similar answers.
It would be interesting and surprising if in each case it is the larger bound
which is closer to the truth.

\subsection{Zeros of $\zeta'(s)$} 
Can one use random matrix theory to   predict 
the horizontal distribution of the real parts 
of the zeros of $\zeta'$? It is known that the Riemann Hypothesis 
is equivalent to the assertion that each non-real zero of $\zeta'(s)$ 
has real part greater than or equal to 1/2. Moreover, if such a zero 
has real part 1/2, then it is also a zero of $\zeta(s)$ (and so a multiple 
zero of $\zeta(s)$). These assertions are the point of departure for
Levinson's work on zeros of the Riemann zeta-function on the critical line.
It would be interesting to know the horizontal 
distribution of these zeros; in particular what proportion of them 
with ordinates between $T$ and $2T$ 
are within $a/\log T$ of the 1/2-line?

 In a similar vein, the Riemann $\xi$-function is real on the 1/2-line 
 and has all of its zeros there (assuming the Riemann Hypothesis).  
 It is an entire function of order 1; because of its functional equation, 
 $\xi(1/2+i\sqrt{z})$ is an entire function of order 1/2.
  It follows that the Riemann Hypothesis 
  implies that all zeros of $\xi'(s)$ are on the 1/2-line. Assuming this 
  to be true, one can ask about the vertical distribution of zeros of $\xi'(s)$, 
  and more generally of $\xi^{(m)}(s)$.  It seems that the zeros of higher 
  derivatives will become more and more regularly spaced; can 
  these distributions be expressed in a simple way using  random matrix theory?

\subsection{Long Mollifiers, Local Integrals, and GUE}
David Farmer \cite{F1}, \cite{F2} has made 
two very interesting conjectures having to do with $\zeta(s)$. 
The first is a conjecture about the mean square of $\zeta(1/2+it)$ 
 times an arbitrarily long mollifier. A mollifier is a  Dirichlet polynomial
  with coefficients equal to the M\"{o}bius $\mu$-function times a smooth function.
  The length of the mollifier is the length of the Dirichlet polynomial. 
   He has also conjectured that
    \begin{equation}
     \frac 1 T\int\sb 0\sp T \frac{\zeta(  s+u)\zeta(  1-s+v)}
{\zeta( s+a) \zeta(  1-s+b)} \,dt \sim  1 + \frac{(u-a)(
v-b)} { (u+v)(a+b)}(1 - T\sp {-(u+v)}).
\end{equation}
       where $a,b,u,v$ are
      complex numbers with positive real part, and $s=1/2+it$.  
 These two conjectures are essentially equivalent and imply certain parts of the GUE conjecture. 
It would be interesting to generalize these and relate them to the full GUE conjecture.

\section{Appendices}

\subsection{  The Classical Groups}
\begin{itemize}
\item The {\bf unitary group} U($N$)  is the group of $N\times N$ matrices $U$ with entries in 
${\Bbb C}$ for which $UU^*=I$ where $U^*$ denotes the conjugate transpose of $U$, i.e. if 
$U=(u_{i,j})$, then $U^*=(\overline{u_{j,i}})$.
\item The {\bf orthogonal group} O($N$)  is the subgroup of U($N$) consisting of matrices with 
real entries.
\item The {\bf special orthogonal group} SO($N$). This is the subgroup of O($N$) consisting of matrices 
with determinant 1.
SO(2$N$) leads to the symmetry type we have called O$^+$ and SO(2$N$+1) leads 
to the symmetry type we call O$^-$.
\item The {\bf symplectic group} USp(2$N$)    is the subgroup of U(2$N$) of matrices $U$ for which 
$UZU^t=Z$ where $U^t$ denotes the transpose of $U$ and 
\begin{equation}Z=\bigg(\begin{array}{cc} 0& I_N\\-I_N&0\end{array}\bigg)\end{equation}
\end{itemize}
 
\subsection{  The Weyl Integration Formula}
The $N\times N$ unitary matrices can be parametrized by their $N$ eigenvalues on the unit circle.
Any configuration of $N$ points on the unit circle corresponds to a conjugacy class of U($N$). 
If $f(A)=f(\theta_1,\dots \theta_N)$ is a symmetric function of $N$ variables,
 then Weyl's formula \cite{W} gives
\begin{equation}
\int_{U(N)}f(A) dA=\frac{1}{N!}\int_{[0,1]^N}f(\theta)\prod_{1\le 
j<k\le N}|e(\theta_j)-e(\theta_k)|^2 ~d\theta_1\dots d\theta_N \end{equation}
where $dA$ is the Haar measure.
Similarly, on $\mbox{Sp}(2N)$ and SO($2N$) we have respectively
\begin{equation}  dA = \frac{2^{N ^2}}{N!}  \prod_{j<k}(\cos \pi\theta_j-\cos\pi \theta_k)^2
 \prod_{j=1}^N \sin ^2 \pi \theta_j\prod_{j=1}^N d\theta_j \thinspace;\end{equation}
\begin{equation}  dA = \frac{2^{(N-1)^2}}{N!} \prod_{i<j}(\cos 
\pi \theta_j-\cos \pi \theta_k)^2
\prod_{j=1}^N d\theta_j  \thinspace.\end{equation}

\subsection{Four Statistics}

Suppose we have a  sequence ${\mathcal T} $ of $N$-tuples of numbers $T_N=\{t_1, t_2,  \dots , t_N\}$ 
where $t_1\le t_2< \dots < t_N$ such that for each set the average spacing $t_{j+1}-t_j$ is 
asymptotically 1.  We write $t_{i,N}$  in place of $t_i$ if we need to indicate that $t_i\in T_N$.

\begin{itemize}
\item The {\bf $n$-level density} of $\mathcal{T}$ is $W(\vec{x})=W(x_1,\dots x_n)$
 means that
\begin{equation}\lim_{N\to \infty}  \sum_{({i_1},\dots {i_n}),i_j\le N
\atop  {i_j}\ne  {i_k}}f(t_{i_1}, 
\dots t_{i_n})
=  \int_{{\Bbb R}^n}f(\vec{x})W(\vec{x}) ~d\vec{x} \thinspace.
\end{equation}
 for a   Schwarz-class 
  $f$. The
 sum 
 is over  $n$-tuples with distinct entries. 
\item The {\bf $j$-th lowest zero density} is $\nu_j(x)$ means that for a test function $f$
\begin{equation}\lim_{N\to \infty}\frac 1 N \sum_{n\le N}f(t_{j,n}) = \int_0^\infty
f(x)\nu_j(x)~dx\thinspace.\end{equation}
\item The {\bf consecutive spacing density} is $\mu(x)$ means that for a test function
$f(x)$ we have
\begin{equation}\lim_{N\to \infty}\frac
 1 N \sum_{i\le N-1}f(t_{i+1}-t_i)=\int_0^\infty f(x)\mu(x)~dx\thinspace.\end{equation}
\item The {\bf $n$-correlation density} is $V(x_1, \dots, x_n)$
means that for   test functions $f$ that are symmetric in all of the variables, depend only 
on the differences of the variables (i.e. 
$ f(x_1+u, \dots , x_n+u)=f(x_1, \dots, x_n)$ for all $u$), and are rapidly decaying on the 
hyperplane $P_n:\{(x_1,\dots, x_n:\sum x_i=0\}$, we have, as $N\to \infty$,
\begin{equation}  \frac 1{N } \sum_{t_1,\dots t_n\in T_N \atop  {i_j}\ne  
{i_k}}f(t_1, \dots t_n)
 \to \int_{P_n}f(\vec{x})V(\vec{x}) ~dx_1\dots ~dx_{n-1}\end{equation}
\end{itemize}
as $N\to \infty$.
The spacing and $n$-correlation densities   are universal, i.e. the 
same for each of O, Sp, and U,
whereas the $n$-level and $j$-th lowest zero densities depend on the symmetry type.

\subsection{ Gaudin's Lemma}
 Associated to each $N\times N$ unitary matrix 
  $A$ are its $N$ eigenvalues  
$\mbox{e}(\theta_j)$ where $0\le \theta_1 \le \dots \le \theta_N\le 1$.
We integrate a function $F(A)$ over U($N$)    by parametrizing the 
group by the $\theta_i$ and using Weyl's formula to convert the integral 
into an $N$-fold integral over the $\theta_i$. 

Often  one wants to integrate with respect to Haar measure over U($N$) 
a function $F(A)=\tilde{f}(A)=\tilde{f}(\theta_1,\dots,\theta_N)$ of $N$ variables that is ``lifted '' from a function $f$ of $n$ variables:
\begin{equation}
\tilde{f}(\theta_1,\dots,\theta_N)=\sum_{(i_1,\dots,i_n) \atop
i_j\ne i_k} f(\theta_{i_1},\dots,\theta_{i_n}) 
\end{equation}
where the sum is over all possible $n$-tuples $(i_1,\dots,i_n)$ 
of distinct integers between 1 and $N$. 
 Gaudin's lemma gives a simplification of this computation 
from an $N$-fold integral to an $n$-fold integral.
The  Haar measure (see section 5.2) at the matrix $A$ can be expressed as
\begin{equation}dA=
 \frac{1}{N!}  \prod_{1\le 
j<k\le N}|e(\theta_j)-e(\theta_k)|^2 ~d\theta_1\dots d\theta_N\thinspace.
\end{equation}
The product here is the square of the absolute value of the  $N\times N$ Vandermonde determinant 
with $j,k$ entry $\mbox{e}(\theta_k)^{j-1}=\mbox{e}((j-1)\theta_k)$ .
It is also the $N\times N$ determinant of the matrix with $j,k$ entry $J_N(\theta_j-\theta_k)$ 
   where
\begin{equation}
J_N(\theta)=\sum_{m=0}^{N-1}\mbox{e}(m\theta)=\mbox{e}((N-1)\theta/2)\frac{\sin \pi N \theta}{\sin \pi \theta}\thinspace.
\end{equation}
Thus,
\begin{equation}
\int_{U(N)}\tilde{f}(A)dA=\int_{[0,1]^N}\tilde{f}(\theta_1,\dots,\theta_N)
\frac{1}{N!}\det_{N\times N}J_N(\theta_j-\theta_k)\prod_{j=1}^Nd\theta_j\thinspace.
\end{equation}
Then, Gaudin's lemma 
asserts the equality
\begin{equation} \int_{ {U}(N)}\tilde{f}(A) dA =\int_{[0,1]^n} f(\theta_1,\dots,\theta_n)
\frac{1}{n!}\det_{n\times n}J_N(\theta_j-\theta_k)\prod_{j=1}^nd\theta_j\thinspace.
 \end{equation}

This principle works for all of the subgroups of U($N$) under consideration here as well.
(See [KS] section 5.1 for a general statement and proof of this important lemma.)
We illustrate by  computing  the $n$-level density function for U($N$). 
Note that 
\begin{equation}
\lim_{N\to \infty} \frac{1}{N}J_N(\theta/N)=\mbox{e}(\theta/2)\frac{\sin \pi \theta}{\pi \theta}
\end{equation}
from which it follows easily that
\begin{equation}
\lim_{N\to \infty} \frac{1}{N^n}  \det_{n\times n} J_N(\theta_j-\theta_k)
= \det_ {n\times n} K_{0}(\theta_1,\dots,\theta_n)
\end{equation} where $K_\epsilon$ is defined  in section 5.5.

Now let $f(\vec{x})=f(x_1,\dots,x_n)$ be a test function.  To compute the $n$-level density (compare
 with section 5.4) we need to evaluate
 \begin{equation}
 \lim_{N\to \infty}   \int_{{U}(N)}\sum_{(i_1,\dots,i_n) \atop
i_j\ne i_k} f(\hat{\theta}_{i_1},\dots,\hat{\theta}_{i_n})\prod_{j<k}|\mbox{e}(\theta_j)-\mbox{e}(\theta_j)|^2~d\theta_1\dots~d\theta_N
\end{equation}
 By Gaudin's lemma and after using the new expression for the Haar measure and
 changing variables $\theta_j\to x_j/N$, the   above is equal to
\begin{eqnarray}
 \lim_{N\to \infty}&&\frac {1}{N^n}\int_{[0,N]^n}f( {\theta}_1,\dots, {\theta}_n)
 \frac{1}{n!}\det_{n\times n}J_N(\theta_j-\theta_k)\prod_{j=1}^nd\theta_j\\
&&= \int_{\Bbb{R}^n}f(x_1,\dots,x_n)\det_{n\times n}{K_0}(x_1,\dots,x_n ) ~dx 
\end{eqnarray}
 so that $W_{ {U},n}(x_1,\dots,x_n)=\det_{n\times n } K_0(x_1,\dots,x_n)$.

\subsection{ Formulas for the Density Functions}
\begin{itemize}
\item
The {\bf $n$--level density} is $W(x_1, \dots, x_n)=\det_{n \times n}
K_\epsilon(x_1,\dots,x_n)$ 
 where $K_\epsilon(x_1,\dots,x_n)$ is the $n\times n$ matrix with entries
\begin{equation}(K_\epsilon(x_1,\dots,x_n))_{i,j}=\frac{\sin \pi(x_i-x_j)}{\pi(x_i-x_j)}+\epsilon \frac{\sin 
\pi(x_i+x_j)}{\pi(x_i+x_j)}
\end{equation}
where $\epsilon=0$ for U; $\epsilon=-1$ for Sp;
 $\epsilon=1$ for O$^+$.
Also, 
\begin{equation}
W_ {O^-,n}(\vec{x}) =\det_{n \times n} (K_{-1}(\vec{x})) +
 \sum_{m=1}^n \delta(x_m) \det_{n-1\times n-1}  ( K^{(m)}_{-1}(\vec{x}))
\end{equation}
where $\delta$ is the  Dirac $\delta$-function
 and the superscript $m$ denotes that the $m$-th row and $m$-th column have been deleted from
  $K_{-1}(\vec{x})$.
\item 
The {\bf lowest zero density} is $\nu_1(x)$ where
\begin{equation}
\nu_1(x)=-\frac{d}{dx}\prod_{j=0}^\infty (1-\lambda_j(x)) \qquad \mbox{U}\thinspace;
\end{equation}
\begin{equation} 
\nu_1(x)= -\frac{d}{dx}\prod_{j=0}^\infty (1-\lambda_{2j+1}(2x)) \qquad  \mbox{Sp}\thinspace;
\end{equation}
\begin{equation}
 \nu_1(x)= -\frac{d}{dx}\prod_{j=0}^\infty (1-\lambda_{2j}(2x)) \qquad \mbox{O}\thinspace,
 \end{equation} 
 where $1\ge \lambda_0(x)\ge \lambda_1(x)\dots $
 are the eigenvalues of
\begin{equation}\int_{-x/2}^{x/2}\frac{\sin 
\pi(t-u)}{\pi(t-u)}f(u)~du=\lambda(x)f(t)\end{equation} 
\item The {\bf consecutive spacing density} is 
\begin{equation}
 \mu(x) 
=\prod_{j=0}^\infty (1-\lambda_j(x))\thinspace.\end{equation}
\item The {\bf $n$--correlation density}, $V(x_1,\dots , x_n)=W_{U,n}(x_1,\dots,x_n)$,

\end{itemize}

\subsection{ The Selberg Integral}

There are many versions of Selberg's integral see \cite{Me}; one is as follows.

If $\Re \alpha>0, \Re \beta >0, \Re \gamma >-\min(\frac 1 n,\frac{\Re \alpha}{n-1},\frac{\Re 
\beta}{n-1})$, then
\begin{eqnarray} 
&&\int_{-1}^1\dots\int_{-1}^1 \prod_{1\le i<j\le 
N}|x_i-x_j|^{2\gamma}\prod_{j=1}^n(1-x_j)^{\alpha-1}(1+x_j)^{\beta-1}dx_j\\
&&\quad=2^{\gamma n(n-1)+n(\alpha+\beta-1)}\prod_{j=0}^{n-1} 
\frac{\Gamma(1+\gamma+j\gamma)\Gamma(\alpha+j\gamma)\Gamma(\beta+j\gamma)}{\Gamma(1+\gamma)\Gamma(
\alpha+\beta+\gamma(n+j-1))}\thinspace.
\end{eqnarray}

\subsection{Moments of Characteristic Polynomials}
\begin{eqnarray}
M_{U,N}(s)&=&\int_{U(N)} |\det(A-I\exp(-i x))|^{2s} 
~dA\\&=&\prod_{j=1}^N\frac{\Gamma(j)\Gamma(j+2s)}{\Gamma(j+s )^2}\thinspace,
\end{eqnarray}
 \begin{eqnarray}
M_{Sp,2N}(s)&=&\int_{Sp(2N)} |\det(A-I)|^{s} 
~dA\\&=&2^{2Ns}\prod_{j=1}^N\frac{\Gamma(1+N+j)\Gamma(1/2+s+j+s)}
{\Gamma(1/2+j)\Gamma(1+s+N+j)}\thinspace,
\end{eqnarray}
\begin{eqnarray}
M_{O,2N}(s)&=&\int_{O(2N)} |\det(A-I )|^{s} 
~dA\\&=&2^{Ns}\prod_{j=1}^N\frac{\Gamma(N+j-1)\Gamma(s+j-1/2)}
{\Gamma(j-1/2)\Gamma(s+j+N-1)}\thinspace.
\end{eqnarray}

\end{document}